\documentclass{amsart}
\usepackage{graphicx}

\newtheorem{proposition}{Proposition}
\newtheorem{theorem}{Theorem}

\theoremstyle{definition}

\newtheorem{example}{Example}[section]
\newtheorem{definition}{Definition}

\newcommand{\myA}{L}
\newcommand{\myB}{S}
\newcommand{\myK}{\mathbf{\Delta Gr}}

\newcommand{\myP}{\mathbf{Topp}}
\newcommand{\myR}{\mathbf{Rings}}
\title{On  ring-valued invariants of topological pairs}
\author{R. M.  Kashaev}
\address{Universit\'e de Gen\`eve,
Section de math\'ematiques,
2-4, rue du Li\`evre,
CP 64,
1211 Gen\`eve 4, Suisse
}
\date{21 January 2007}
\email{Rinat.Kashaev@math.unige.ch}
\begin{document}
\maketitle
\begin{abstract}
Based on the notion of a $\Delta$-group(oid), ring-valued invariants of pairs of topological spaces can be defined in intrinsic topological terms.
\end{abstract}
\section{Introduction}

Let $\myP$ be the category of quadruples  $(X,A,L,\delta)$,  where $A\subset X$ is a pair of topological spaces, $L$ is a subset  of the set $\pi_1(X,A)=[I,\partial I;X,A]$ of homotopy  classes of paths $\gamma \colon (I,\partial I)\to (X,A)$ (see, for example, \cite{Hatcher} for definitions)  such that if  $[\gamma]\in L$ then $[\bar{\gamma}]\in L$, where
$\bar{\gamma}(t)=\gamma(1-t)$, and $\delta=\{\delta_0,\delta_1\}$ is a pair of inclusions 
\[
\delta_i\colon\pi_1(X,A) \to V_\delta \subset A,\quad i\in\{0,1\}
\]
 satisfying the following conditions:
\begin{itemize}
\item[(i)] for any class $[\gamma]\in \pi_1(X,A)$ there exists a representative path $\alpha$, such that $\alpha(i)=\delta_i([\gamma])$, $i\in\{0,1\}$;
\item[(ii)]  $\delta_0([\gamma])=\delta_1([\gamma])$ if and only if $[\gamma]=[\bar{\gamma}]$.
\end{itemize}
A morphism $f\colon (X,A,L,\delta)\to (Y,B,M,\rho)$ in $\myP$ is a continuous map of triples of topological spaces  $f\colon(X,A,V_\delta)\to(Y,B,V_\rho)$ such that $f_*(L)\subset M$ and $f\circ\delta_i=\rho_i\circ f_*$.
The main purpose of this paper is to construct a non-trivial functor $r\colon\myP\to\myR$.

The construction of this functor goes as follows.  First,  we define a category $\myK$ of  $\Delta$-groupoids (and, in particular,   $\Delta$-groups). We observe that the group of invertible elements of any ring admits a canonical structure of a   $\Delta$-group. Then, we show that there is a functorial way of associating   $\Delta$-groupoids to objects of $\myP$ (Theorem~\ref{myth:2}). Finally,  we define the functor  $r$ through a universal object construction.
 
\section{  $\Delta$-groupoids and  topological pairs}\label{mysec:1}
\subsection{$\Delta$-group(oid)s}
Let $G$ be a groupoid and $H$ its subset. We say that a pair of elements $(x,y)\in H\times H$ is $H$-\emph{composable} if it is composable in $G$ and $xy\in H$. 
\begin{definition}
A \emph{  $\Delta$-groupoid} is a triple $(G,H,k)$, where $G$ is a groupoid, $H$ a subset of $G$ closed under the inversion map $i\colon H\ni x\mapsto x^{-1}\in H$, and $k\colon H\to H$ an involution such that 
for any $H$-composable pair $(x,y)$ the pairs $(k(xy),ik(y))$ and $(k(x),iki(y))$ are also $H$-composable,  and the following identity holds true: 
\begin{equation}\label{myeq:1}
k(xy)ik(y)=k(k(x)iki(y)).
\end{equation}
In the  case when $G$ is a group, the    $\Delta$-groupoid $(G,H,k)$ is called a \emph{  $\Delta$-group}.
\end{definition}
A morphism $f\colon (G,H,k)\to(G',H',k')$ between two $\Delta$-groupoids is a morphism of groupoids $f\colon G\to G'$ such that $f(H)\subset H'$ and $k'(f(x))=f(k(x))$, $\forall x\in H$. In this way we come to the category $\myK$ of   $\Delta$-groupoids. In what follows, sometimes we shall write $G$ instead of $(G,H,k)$ if the structure of a   $\Delta$-group(oid) is clear from the context.
 
It is immediate to see that in any   $\Delta$-groupoid $(G,H,k)$ we have the identity $iki=kik$. Indeed, applying  identity
 \eqref{myeq:1} to the $H$-composable pair $(k(x),iki(y))$ in its right hand side, we obtain
 \[
 k(xy)ik(y)=k(k(x)iki(y))=k(xy)kiki(y)
 \]
 which implies that $ik(y)=kiki(y)$.
 \begin{example}\label{myex:-1}
 Any group(oid) $G$ is trivially a $\Delta$-group(oid) as the triple 
  $(G,\emptyset,1_\emptyset)$.
  \end{example}
 \begin{example}\label{myex:0}
 Let $X$ be a set. Define  a groupoid
 \[
 X^3=\{(a,b,c)\vert\ a,b,c\in X\},
 \]
 where a pair $((a,b,c),(e,f,g))$ is composable if and only if $a=e$ and $c=f$ with the product
 $(a,b,c)(a,c,g)=(a,b,g)$. The units are of the form $(a,b,b)$, and the inversion map $(a,b,c)^{-1}=(a,c,b)$. For the set
 \[
 X^{3*}=\{(a,b,c)\vert\ a\ne b\ne c\ne a\}\subset X^3,
 \]
 we define  a map $k\colon X^{3*}\ni (a,b,c)\mapsto(c,b,a)\in X^{3*}$. Then, the triple $(X^3,X^{3*},k)$ is a   $\Delta$-groupoid.
 \end{example}

\begin{example}\label{myex:1}
Let $R$ be a ring, and an involution $k\colon R\to R$ be defined by the formula $k(x)=1-x$. Then, the triple $(R^*, k(R^*)\cap R^*,k)$, where $R^*$ is the group of invertible elements of $R$, is a   $\Delta$-group. Thus, we observe that the group of invertible elements of any ring admits a canonical structure of a   $\Delta$-group.
\end{example}
\begin{example}\label{myex:2}
For a ring $R$ let the set $R^*\times R$ be given a group structure with the multiplication $(x,y)(u,v)=(xu,xv+y)$, the unit element $(1,0)$, and the inversion map $(x,y)^{-1}=(x^{-1},-x^{-1}y)$. Then, the triple $(R^*\times R,R^*\times R^*,k)$, where $k(x,y)=(y,x)$, is a   $\Delta$-group. 
\end{example}
For a given ring $R$, there is a canonical morphism
of   $\Delta$-groups from the associated   $\Delta$-group of Example~\ref{myex:1} to that of Example~\ref{myex:2} given by the formula $x\mapsto (x,1-x)$.
\begin{example}
Let $(G,G_\pm,\theta)$ be a symmetrically factorized group of \cite{KR}. That means that $G$ is a group with two isomorphic subgroups $G_\pm $ conjugated to each other by an involutive element $\theta\in G$, and the restriction of the multiplication  map $m\colon G_+\times G_-\to G_+G_-\subset G$ is a set-theoretical bijection, whose inverse is called the factorization map $G_+G_-\ni g\mapsto (g_+,g_-)\in G_+\times G_-$. In this case, the triple $(G_+,G_+\cap G_-G_+\theta\cap\theta G_+G_- ,k)$, where $k(x)=(\theta x^{-1})_+^{-1}$, is a   $\Delta$-group. Moreover, from the work \cite{KR} it follows, that any $\Delta$-group corresponds to a symmetrically factorized group.
\end{example}
\subsection{The  functor  $g\colon \myP\to \myK$}
Let $X=(X,A,L,\delta)\in\mathrm{Ob}(\myP)$. An element of $\myA$ will be called a \emph{long arc}.  Any long arc can be identified with an element  of the fundamental groupoid $\pi_1(X,V_\delta)=[I,\partial I;X,V_\delta]$ through the natural identifications of the object maps: $[\gamma](i)=\delta_i([\gamma])$.

Let 
\[
\pi_1(A\subset X,V_\delta)=i_{A*}(\pi_1(A,V_\delta))\subset \pi_1(X,V_\delta)
\]
be the subgroupoid of path classes in $A$, where  $i_A\colon A\to X$ is the inclusion map. We define a subset $\myB_X\subset \pi_1(A\subset X,V_\delta)$ of  \emph{short arcs} as follows. A class $x\in \myB_X$ if and only if there exist $y,z\in \pi_1(A\subset X,V_\delta)$ and long arcs 
 $\alpha,\beta,\gamma$ such that the triples $(\alpha, x, \beta)$ and $(y,\gamma,z)$ are composable in $\pi_1(X,V_\delta)$, and satisfy the equation 
\begin{equation}\label{myeq:2}
\alpha x\beta=y\gamma z.
\end{equation}
It is clear, that $y$ and $z$ are also short arcs. From the definition of a long arc, it follows that for a given $x\in\myB_X$, the triple $(\alpha,\beta,\gamma)$, satisfying equation~\eqref{myeq:2}, is unique. Thus, we can define three maps\footnote{The first two of these maps, $p_0$ and $p_1$, are, in fact, restrictions of the maps $p_i\colon\pi_1(A\subset X,V_\delta)\to\pi_1(X,A)$ defined through the equations $p_i(x)=\delta_{1-i}^{-1}( x(i))$.} $p_i\colon  \myB_X\to \myA$, $i\in\{0,1,2\}$, by the equations $p_0(x)=\alpha$, $p_1(x)=\beta$, $p_2(x)=\gamma$. 

On the other hand, for a given $x\in\myB_X$, the pair $(y,z)$, satisfying equation~\eqref{myeq:2}, is not necessarily unique. For this reason,  we define a minimal equivalence relation $\sim$ on the set $\myB_X$ which implies uniqueness of the pair $([y],[z])$ of equivalence classes, i.e. if  for  a given $x\in \myB_X$, there are two pairs $(y,z)$ and $(y',z')$  of short arcs, satisfying equation ~\eqref{myeq:2}, then $y\sim y'$ and $z\sim z'$.

Let $G_X$ be a groupoid obtained from $\pi_1(A\subset X,V_\delta)$ by imposing additional relations of the form $x=y$, where $x\sim y\in \myB_X$.  The quotient set $H_X=\myB_X/\sim$ is naturally identified as a subset of $G_X$ closed under the inversion map $i$. We define another map $k=k_X\colon H_X\to H_X$ by assigning $k([x])=[z]$, where $x$ and $z$ satisfy equation~\eqref{myeq:2}. Clearly, $k$ is an involution. Besides,  by the symmetry of the equation~\eqref{myeq:2}, one easily sees that $ki([x])=i([y])$, $k([y])=i([z])$. Altogether, these imply that $kiki([x])=ik([x])$, i.e. $iki=kik$. Thus, equation~\eqref{myeq:2}  leads to an equation of the form
\begin{equation}\label{myeq:3}
p_0(  x)   x p_1(  x)=j(  x)p_2(  x) k(  x),\quad j=iki=kik, \quad \forall x\in H_X.
\end{equation}
\begin{theorem}\label{myth:2}
Let  $(X,A,L,\delta)\in\mathrm{Ob}(\myP)$. Then, the triple 
\[
g(X,A,L,\delta)=(G_X,H_X,k_X)
\]
 is a   $\Delta$-groupoid. Moreover, the map $g\colon \myP\to\myK$ is  a covariant functor.
\end{theorem}
\begin{proof}
To simplify notation, below we shall not distinguish between the short arcs and the corresponding equivalence classes in $H_X$. Let $(x,y)$ be an $H_X$-composable pair. Then, $xy(1)=y(1)$, and thus $p_1(xy)=p_1(y)$. From the definition of the map $k$, it follows that for any $z\in H_X$, if $\alpha=p_1(z)$, then $p_1(k(z))=\alpha^{-1}$. This implies that $p_1(k(xy))=p_1(k(y))$, thus, 
the pair $(k(xy),ik(y))$ is composable.
Similarly, the two other pairs $(ki(xy),iki(x))$ and $(k(x),iki(y))$ are also composable. Denote the corresponding compositions as $u=k(xy)ik(y)$, $v=ki(xy)iki(x)$, and $w=k(x)iki(y)$. Now, by using consecutively equation~\eqref{myeq:3}, we have the following sequence of  equalities
\begin{multline*}
p_0(u) u  p_1(u)= p_0k(xy) k(xy) ik(y)  p_1ik(y)\\
=p_0k(xy)  k(xy)  p_1k(xy)  p_0ik(y) ik(y)  p_1ik(y)\\
= jk(xy)  p_2k(xy)   xy jik(y)  p_2ik(y)  j(y) =ki(xy)  p_0(x)   x  p_1(x)  j(y)\\
= ki(xy)j(x)  p_2(x)   k(x)j(y)=v  p_2(u)  w
\end{multline*}
which implies that $u,v,w\in H_X$, and $w=k(u)$. The latter  is equivalent to equation~\eqref{myeq:1}.

Functoriality of $g$ easily follows from the fact that a morphism 
\[
f\colon (X,A,L,\delta)\to (Y,B,M,\rho)
\]
 in $\myP$ induces a morphism of groupoids
\[
f_*\colon \pi_1(A\subset X,V_\delta)\to \pi_1(B\subset Y,V_\rho)
\]
 which sends $\myB_X$ to $\myB_Y$, and equivalent short arcs in $\myB_X$ to equivalent short arcs
 in $\myB_Y$.
\end{proof}
\begin{proposition}
For any $\delta,\rho$, the $\Delta$-groupoids $g(X,A,L,\delta)$ and $g(X,A,L,\rho)$ are isomorphic.
\end{proposition}
\begin{proof}
We have a well defined bijection \( f=\rho_0\circ\delta_0^{-1}=\rho_1\circ\delta_1^{-1}\colon V_\delta\to V_\rho\) such that for any $x\in V_\delta$ there exists a  path $\lambda_x\colon (I,0,1)\to (A,x,f(x))$. For a family $\Lambda$ of such paths we obtain the following groupoid isomorphism
\[
f_\Lambda \colon\pi_1(A\subset X,V_\delta)\to\pi_1(A\subset X,V_\rho),\quad
f_\Lambda (\alpha)=\bar{\lambda}_{\alpha(0)}\cdot\alpha\cdot{\lambda}_{\alpha(1)},
\]
which, clearly, is consistent with the equivalence relations on the corresponding sets of short arcs.
\end{proof}
\subsection{The functor $r\colon \myP\to\myR$}

Now, we define the ring $r(X,A,L,\delta)$ as the quotient ring $\mathbb{Z}G_X/\mathcal{I}_X$, where the  two-sided ideal $\mathcal{I}_X$ is generated by the set of elements $\{k(x)+x-1\vert\ x\in H_X\}$. By construction, we have a   $\Delta$-groupoid morphism $h\colon g(X,A,L,\delta)\to r(X,A,L,\delta)^*$ (see Example~\ref{myex:1} for the canonical $\Delta$-group structure on the group of invertible elements of an arbitrary ring) which is universal in the sense that for any ring $R$ and any    $\Delta$-groupoid morphism $s\colon  g(X,A,L,\delta)\to R^*$ there is a unique ring homomorphism $f\colon r(X,A,L,\delta)\to R$ such that $s=f\circ h$. 

\subsection{Examples}
In the following examples the subset $L\subset\pi_1(X,A)$ will be chosen as the maximal set of nontrivial classes \footnote{A homotopy class of maps of pairs of topological spaces is \emph{nontrivial} if it does not contain a constant map \cite{Hatcher}.}.
\begin{example}
If $(X,\emptyset,L,\delta)\in \mathrm{Ob}(\myP)$, then, obviously, $g(X,\emptyset,L,\delta)=\emptyset$, 
and thus $r(X,\emptyset,L,\delta)=0$.
\end{example}
\begin{example}
Let $(X,X,L,\delta)\in \mathrm{Ob}(\myP)$ with path-connected $X$. Then, obviously, $\pi_1(X,X)=\{[c_{x_0}]\}$, where $c_{x_0}$ is the constant path  at $x_0$.  This implies that  $V_\delta=\{x_0\}$ with $\delta_i([c_{x_0}])=x_0$.
In this case $\myB_X=\emptyset$, $g(X,X,L,\delta)=(\pi_1(X,x_0),\emptyset,1_{\emptyset})$  and thus $r(X,X,L,\delta)=\mathbb{Z}\pi_1(X,x_0)$.
\end{example}

\begin{example}
Let $(X,A(n),L,\delta)\in \mathrm{Ob}(\myP)$ with simply connected $X$, and $A(n)$ consisting of $n$ simply connected path-components. With the notation  $\underline{n}=\{i\in\mathbb{Z}\vert\ 1\le i\le n\}$, we have the bijections $\pi_1(X,A(n))\simeq \underline{n}^2$,  $g(X,A(n),\delta)\simeq \underline{n}^3$ (see Example~\ref{myex:0} for the definition of this $\Delta$-groupoid). 
For the first three $n$ it is straightforward to verify the following ring isomorphisms: $r(X,A(1),L,\delta)\simeq\mathbb{Z}$,  $r(X,A(2),L,\delta)\simeq\mathbb{Z}F_{\{x_0,x_1\}}$, $r(X,A(3),L,\delta)=\{x_0-1\}^{-1}\mathbb{Z}F_{\{x_0,x_1,x_2,x_3\}}$. Here, we use the following notation: $F_S$ is the free group generated by the set $S$; $S^{-1}R$ is the ring $R$ localized with respect to the multiplicative system generated by the set $S$.

\end{example}

\subsection*{Acknowledgments} I would like to thank V. Florens, J.-C. Hausmann, C. Weber, and T. Vust for helpful discussions. The work is supported in part by the Swiss National Science
  Foundation.

\end{document}